\documentclass[11pt,palatino]{article}

\usepackage[numbers,sort&compress,square]{natbib}

\usepackage{amsthm, amsmath, mathtools, mathrsfs}
\usepackage{algorithm}
\usepackage{nccmath}
\usepackage{comment}

\usepackage{amsfonts,amsmath,amssymb,bm,url}

\usepackage{graphicx,epsfig,subfigure}

\usepackage[utf8]{inputenc} 
\usepackage[T1]{fontenc}    
\usepackage{hyperref}       
\usepackage{url}            
\usepackage{graphicx}
\usepackage{booktabs}       
\usepackage{amsfonts}       
\usepackage{nicefrac}       
\usepackage{microtype}      
\usepackage{xcolor}         
\usepackage{comment}
\usepackage{amsmath}
\usepackage{amsthm}
\usepackage[utf8]{inputenc}
\usepackage[english]{babel}
\usepackage{authblk}

\usepackage{comment}

\newcommand*\FF{{\bf F}}

\newcommand*\AAA{{\bf A}}
\newcommand*\BBB{{\bf B}}
\newcommand*\CCC{{\bf C}}

\newcommand{\mcN}{\mathcal{N}}
\newcommand{\mcD}{\mathcal{D}}
\newcommand{\pmcD}{\partial\mathcal{D}}

\newcommand{\mcO}{\mathcal{O}}
\newcommand{\mcL}{\mathcal{L}}

\newcommand{\mcC}{\mathcal{C}}
\newcommand{\tsinh}{\text{sinh}}
\newcommand{\tcosh}{\text{cosh}}
\newcommand{\hmcNe}{\hat{\mathcal{N}}_\epsilon}
\newcommand{\mcNe}{{\mathcal{N}}_\epsilon}

\begin{document}

\title{Error Estimation and Correction from within Neural Network Differential Equation Solvers}

\author[1, 2]{Akshunna S. Dogra}
\affil[1] {\small Department of Mathematics, Imperial College London}
\affil[2] {\small EPSRC CDT in Mathematics of Random Systems: Analysis, Modelling and Simulation}
\date{}

\maketitle

\begin{abstract}
Neural Network Differential Equation (NN DE) solvers have surged in popularity due to a combination of factors: computational advances making their optimization more tractable, their capacity to handle high dimensional problems, easy interpret-ability of their models, etc. However, almost all NN DE solvers suffer from a fundamental limitation: they are trained using loss functions that depend only implicitly on the error associated with the estimate. As such, validation and error analysis of solution estimates requires knowledge of the true solution. Indeed, if the true solution is unknown, we are often reduced to simply hoping that a "\textit{low enough}" loss implies "\textit{small enough}" errors, since explicit relationships between the two are not available/well defined. In this work, we describe a general strategy for efficiently constructing error estimates and corrections for Neural Network Differential Equation solvers. Our methods do not require advance knowledge of the true solutions and obtain explicit relationships between loss functions and the error associated with solution estimates. In turn, these explicit relationships directly allow us to estimate and correct for the errors.
\end{abstract}

\section{Basic assumptions and groundwork}
\subsection{Setup}
Let $\mcD \subseteq \mathbb{R}^d$ denote a closed, bounded, path connected domain of interest for some differential equation. Let $\pmcD \subseteq \mcD$ be the portion of the domain over which some constraint conditions exist (often boundary conditions and/or empirical data). Assume that our chosen differential equation (DE), when given unique constraints over $\pmcD$, admits a unique, analytic solution $\Phi:\mcD \to \mathbb{R}^D$.
We denote the DE as:
\begin{equation}
    \label{Lin_nonlin_ind_eq}
    \FF[\Phi] = \AAA[\Phi] + \BBB[\Phi] + \CCC = 0
\end{equation}
where $\AAA$ represents the term(s) on which $\Phi$ depends linearly, $\BBB$ the represents the term(s) on which $\Phi$ depends nonlinearly, and $\CCC$ are the term(s) independent of $\Phi$. This additive decomposition into linear, nonlinear, and independent terms is always possible: for a linear DE, $\BBB \equiv 0 \equiv \CCC$. Just like $\Phi$, if the operators $\AAA, \BBB, \CCC$ make use of any functions over $\mcD$, we assume those functions are analytic over $\mcD$ as well.

Let us assume we wish to construct a feed forward Neural Network (NN) based approximation $\mcN: \mcD \to \mathbb{R}^D$. We assume the NN uses analytic activation functions as well: as such, $\mcN \in \mcC^{\infty}$. Let $W$ be its width (neurons per hidden layer) and $D_{\mcN}$ be its depth (number of hidden layers). Let ${\bf w}{}\equiv \{b^1_1,w^1_{11}, w^1_{12}, ....,b^2_1,w^2_{11},w^2_{12},....\}\in\mathcal{M}\subseteq\mathbb{R}^M$ be the $M$ weights/biases of the NN. 

\newpage
\subsection{Existence:} Feedforward NNs with ReLU activation and $W\geq d+4$ can arbitrarily approximate any Lebesgue integrable function $\Phi: \mathcal{D} \to \mathbb{R}^D$ on a $L^1$ norm, provided that $\mcD$ is some compact subset of $\mathbb{R}^{d}$ and enough depth $D_{\mcN}$ is provided to the NN \cite{Lu17}. NNs with ReLU activation and $D_{\mcN} \geq log_2(d+1)$ can arbitrarily approximate any Lebesgue integrable function $\Phi:\mcD \to \mathbb{R}^D$ on a $L^1$ norm, provided that enough width $W$ is provided to the NN \cite{Arora18}. Hence, there exist $\bf w$ such that $\mathcal{N}$ can arbitrarily approximate any $\Phi \in \mathcal{C}^k$ over $\mcD$, given large enough $(W,D_{\mcN})$.

The jump from showing the theoretically possible to quantifying the needs of the practical happens by leveraging the minimum regularity expected from $\Phi$. Recall $\mathcal{W}^{k,\infty}(\mcD)$ is the Sobolev space of order $k$, on compact $\mcD \subseteq\mathbb{R}^{d}$, on the $L^{\infty}$ norm.

\subsection{Complexity:} Let us assume that $\Phi\in \mathcal{W}^{k,\infty}(\mcD)$. There exists some $\mathcal{N}$ with $(W,D_{\mcN})$ that can $\varepsilon$ - approximate $\Phi:{\mcD}\to \mathbb{R}$, if $W=d+1, \text{ }D_{\mcN}\text{ } = \mathcal{O}({\text{diam}}({\mcD})/\omega^{-1}_f(\varepsilon))^{d}$ \cite{han17}, where
\[
    \omega_f^{-1}(\varepsilon) = \text{sup}\{\delta:\omega_f(\delta)\leq\varepsilon\},  \text{ }\text{ }\text{ }\text{ }\delta=|{\bf x}_1-{\bf{x}}_2|,\text{ }\text{ }\text{ }{\bf x}_1,{\bf{x}}_2\in \mcD.
\]
There also exist $\mathcal{N}$ that $\varepsilon$ - approximate $\Phi:{\mcD}\to \mathbb{R}$, if $M=\mathcal{O}{(ln(\frac{1}{\varepsilon}){\varepsilon}}^{-{d}/{k}}), \text{ }D_{\mcN} = \mathcal{O}(ln(\frac{1}{\varepsilon}))$ \cite{yar17}. 

\subsection{Optimization}
The mean loss function $\mcL$ to train $\mcN$ is usually: (over-line denotes mean over the sampled points)
\begin{equation}
\label{gen_loss_DE}
\begin{aligned}
&\mcL = \overline{||\FF[\mcN(x_1)]||^p} + \overline{||\mcN(x_2) - \Phi(x_2)||^p} \textbf{ }\textbf{ }\textbf{ }\textbf{ }\textbf{ }\textbf{ }\textbf{ }\textbf{ } \textbf{ }\textbf{ }\textbf{ }\textbf{ }\textbf{ }\textbf{ }\textbf{ }\textbf{ }\textbf{ }\textbf{ }\textbf{ }\textbf{ }\textbf{ }\textbf{ }\textbf{ }\textbf{ } x_1  \in \mcD, x_2 \in \pmcD\\
&\textbf{ }\textbf{ } = \overline{||\AAA[\mcN(x_1)] + \BBB[\mcN(x_1)] + \CCC||^{p}} +  \overline{||\mcN(x_2) - \Phi(x_2)||^{p}}\\
\end{aligned}
\end{equation}
Variants of Stochastic Gradient Descent (SGD) are \textit{almost always} capable of \textit{eventually} reaching adequately small loss/error values for large enough $(W,D_\mcN)$ under such loss functions \cite{lecun14,lecun15,Bottou99,Kiwiel01}.  

Clearly, if $\mcL = 0$, the uniqueness of $\Phi$ implies $\mcN = \Phi$. Thus, the idea behind training is that if/as $\mcL \to 0, \mcN \to \Phi$ over $\mcD$. 

\subsection*{Comments on a modified approach}
Sometimes, the constraint conditions may be enforced by using the following parametrization as the approximation to the solution $\Phi$ at some arbitrary $x_1 \in \mcD$:
\begin{equation}
    \mcN(x_1) = \Phi(x_2) + \text{dist}(x_1, x_2)\mcN_0
\end{equation}
where $x_2 \in \pmcD$ is some corresponding point in $\pmcD$ for $x_1$ (an intuitive choice is the nearest point in $\pmcD$). Here, $\mcN_0$ represents the core of the NN based approximation and $\text{dist}(x_1, x_2)$ is some metric that enforces that $\mcN$ is always exact over $\pmcD$, while allowing the NN flexibility to \textit{learn} $\Phi$ elsewhere (ex: $1 - e^{-||x_1 - x_2||}$, as used in \cite{dog20,dogred20_2,dogred20, mar20}). This parametrization can eliminate the need for the second term in Eq. \ref{gen_loss_DE}. Our work is applicable in such scenarios too.

\newpage
\section{The Problems}
Define $\mcN_\epsilon = \Phi - \mcN$. Some simple observations underpin the entirety of this work:

\begin{itemize}
    \item $\FF[\mcN]$ in Eq. \ref{gen_loss_DE} is not explicit in $\mcNe$: there is an ambiguous relationship between $\mcL$ and $\mcNe$ over $\mcD - \pmcD$, where $\Phi$ is unknown. We don't know \textbf{if/how/when} $\mcNe \to 0$ as $\mcL \to 0$.
    \item Thanks to that ambiguity, $\mcNe$ associated with any NN based estimate $\mcN$ is only estimable over $\mcD$, if $\Phi$ is already known over $\mcD$, thus rendering the approximation $\mcN$ superfluous.
    \item Since finding the $\mcN$ such that $\mcL = 0$ is extremely unlikely, optimization resources are spent on minor refinements on $\mcN$, once the NN settles around a region of local minima of $\mcL$.
\end{itemize}
The first problem can be immensely intractable in nature and has to be handled separately for different kinds of DEs. Formal discussion on this aspect of NN DE solvers in some kinds of problems can be found in \cite{shin_karni_20, dogred20_2, mar20, Sirig18, zimmer20, salimova20}. Fortunately, the scope of this work concerns NN DE solvers, where it is assumed that arbitrarily adequate approximations are obtainable in a practical sense: our strong assumptions on $\Phi,\mcD$ allow us to sidestep the entire issue for now.

The second observation has significant implications: since the true accuracy of the classic NN DE solver is \textit{seemingly} not verifiable from within, any approximations it produces are unreliable, even if we may intuitively be able to estimate what level of error a particular value of $\mcL$ corresponds to. \textit{A priori}, we only know the explicit relationship between $\mcNe$ and $\mcL$ over $\pmcD$. There is a scarcity of work on constructing methods to explicitly estimate $\mcNe$ over $\mcD$ using information available to the NN DE solver: usual error comparisons rely on actual knowledge of $\Phi$ over $\mcD$. 

In the context of this work, the third observation leads to a useful application: we know that the uniqueness of $\Phi$ as a solution makes it pretty unlikely $\mcN({\bf w})$ will give us the exact solution once some minima of optimization is reached. Any further training cannot provide meaningful gains in accuracy. However, even a moderately good estimate $\mcNe$, if it could be found using available information, should be immediately useful.

From the discussion above, it is clear that an explicit relationship between $\mcL$ and $\mcNe$ would imply many benefits. We show that is a very achievable goal under the established assumptions.

\section{The Resolution}
Define $\mcN_\epsilon = \Phi - \mcN$. We provide resolution to the problems above by noticing the following:
\begin{equation}
\label{perturb_trick}
\FF[\mcN] = \AAA[\mcN] + \BBB[\mcN] + \CCC = \AAA[\Phi] - \AAA[\mcNe] + \BBB[\Phi] - \BBB'[\mcN, \mcNe] + \CCC = - \AAA[\mcNe] - \BBB'[\mcN, \mcNe]
\end{equation}
since $\AAA[\Phi] + \BBB[\Phi] + \CCC = 0$ by problem setup, and $\BBB'[\mcN,\mcNe] = \BBB(\Phi) - \BBB(\mcN)$ is a new term in the equation obtained using \textit{some} sort of transformation. 

The new terms $\AAA[\mcNe], \BBB'[\mcN, \mcNe]$ are the crux of this work: questions regarding their existence and uniqueness have strong implications on whether internal error analysis and correction by NN DE solvers is possible or not. The linear behavior inherent to $\AAA$ makes $\AAA[\mcNe]$ relatively easy to handle, so we will focus on $\BBB'$ hereon. Do note that $\BBB'$ has to be obtained in such a way that terms dependent on $\Phi$ disappear from Eq. \ref{perturb_trick}: that step is what allows us to later on analyse errors without depending upon exact knowledge of $\Phi$.

As a motivating example, consider the cases when the nonlinearities in $\FF{}$ are due to the presence of degree 2 or higher polynomials. To exemplify more explicitly, let us assume that $\BBB$ is quadratic in $\Phi$, i.e., $\BBB[\Phi] = \Phi \cdot \Phi$. We use:
\[\Phi \cdot \Phi = \mcN \cdot \mcN + \mcNe \cdot \mcNe + \mcN \cdot \mcNe + \mcNe \cdot \mcN
 \implies
 \Phi \cdot \Phi - \mcN \cdot \mcN = (\mcNe \cdot \mcNe + \mcN \cdot \mcNe + \mcNe \cdot \mcN)
\]
Thus, we get that
\[\BBB'[\mcN, \mcNe] = (\mcNe \cdot \mcNe + \mcN \cdot \mcNe + \mcNe \cdot \mcN)
\]

In general, it is often possible to use a perturbation trick to obtain $\BBB'[\mcN,\mcNe]$ as follows:
\[\BBB[\Phi] = \BBB[\mcN + \mcNe] = \BBB[\mcN] + \Big(\BBB_1[\mcN, \mcNe] + \BBB_2[\mcN, \mcNe] + .... \Big)
= \BBB[\mcN] + \BBB'[\mcN, \mcNe]
\]

\[\implies \BBB[\mcN] = \BBB[\Phi] - \BBB'[\mcN, \mcNe]
\]

For example, one commonly possible form, since $\Phi$ is analytic, for $\BBB'[\mcN,\mcNe]$ is:
\begin{equation}
    \BBB'[\mcN,\mcNe] \equiv \mathcal{T}_1(\mcN) \mcNe + \frac{\mcNe^{\dagger}\mathcal{T}_2(\mcN)\mcNe}{2!} + ...
\end{equation}
when obtained using a Taylor expansion, where $\mathcal{T}_i(\mcN)$ are some $i^{th}$ order derivatives. 

In general, the transformation technique itself is not important, as long as it is accurate (we will give an example later that uses two different tricks to obtain the same result). The general idea is that we are trying to replace the $\BBB[\mcN]$ term in the loss function with $\BBB[\Phi]$, \textbf{plus} some other things, so that we may zero out the equation term evaluated at $\Phi$ (see Eq. \ref{perturb_trick}).

With the above, we may rectify one major issue with Eq. \ref{gen_loss_DE}: the loss can now indeed be explicitly related to the approximation error $\mcNe$ over $\mcD$ as :
\begin{equation}
    \label{reframed_loss}
    \mcL = \overline{||\AAA[\mcNe(x_1)] + \BBB'[\mcN(x_1), \mcNe(x_1)]||^{p}} + \overline{||\mcNe(x_2)||^{p}} \textbf{ }\textbf{ }\textbf{ }\textbf{ }\textbf{ }\textbf{ }\textbf{ }\textbf{ }\textbf{ }\textbf{ }\textbf{ }\textbf{ }\textbf{ }\textbf{ } x_1  \in \mcD, x_2 \in \pmcD
\end{equation}

In the cases we consider in this work (which have a wide class of applications in physics, engineering, macro-molecule dynamics, etc), the $\BBB' [\textbf{ }]$ exists and is well behaved over the entire domain (since we assumed $\Phi$ to be analytic over $\mathcal{D}$). The extension to the cases where $\BBB'[\textbf{ }]$ is not well behaved, but only over a finite set of points in $\mcD$ is also trivial, and thus our work shall be applicable there as well. Settings beyond that level of complexity will require a more customized approach that we aim to research in later work. 

Note that Eqn. \ref{perturb_trick} gives us the following new differential equation true over the entirety of $\mcD$:

\begin{equation}
    \label{rearragned_error_eq}
    \AAA[\mcNe] + \BBB'[\mcN, \mcNe] + \FF[\mcN] = 0
\end{equation}

We immediately hit upon a method to estimate $\mcNe$: we use another NN, $\hmcNe$, in the manner we used $\mcN$ with Eq. \ref{gen_loss_DE}, to obtain an approximation for $\mcNe$. More precisely, once the optimization of $\mcN$ has finished, define:
\begin{equation}
\label{Err_crrct_gen}
\mcL_2 = \overline{||\AAA[\hmcNe(x_1)] + \BBB'[\mcN(x_1), \hmcNe(x_1)] + \FF[\mcN(x_1)]||^{p}} + \overline{||\hmcNe(x_2) - \mcNe(x_2)||^p}\textbf{ }\textbf{ }\textbf{ }\textbf{ }\textbf{ } x_1  \in \mcD, x_2 \in \pmcD
\end{equation}

The uniqueness of $\Phi$ implies that if $\mcN$ is fixed, then $\mcNe$ is a unique solution to the DE given by Eq. \ref{rearragned_error_eq}. Therefore, just as $\mcL = 0 \implies \hat{\mcN} = \Phi$, $\mcL_2 = 0 \implies \hmcNe = \mcNe$.

A powerful application is also glimpsed if Taylor expansions are used. We may write:
\[\AAA[\mcNe] + \left[\mathcal{T}_1(\mcN) \mcNe + \frac{\mcNe^{\dagger}\mathcal{T}_2(\mcN)\mcNe}{2!} + ...\right] + \FF[\mcN] = 0
\]

Notice that if we can throw away $\mcO(||\mcNe||^2)$ or higher terms as inconsequential (justifiable during the end of training, if $||\mcNe|| \ll ||\mcN||$), Eq. \ref{rearragned_error_eq} becomes:
\[\AAA[\mcNe] + \mathcal{T}_1(\mcN) \mcNe  + \FF[\mcN] = 0
\]
which roughly implies $\overline{||\mcNe||} \sim \overline{||\FF[\mcN]||}$, something that empirical observations match. We are using an extra assumption here that in the latter stages of optimization, the behavior of $\mathcal{T}_1(\mcN)$ does not vary as much, since the solution estimate $\mcN$ starts settling into a fixed state as changes in $\bf{w}$ diminish.

\subsection{A Demonstration}
Let us exemplify the above discussion with an example DE that contains terms with linear, nonlinear, and no dependence on the input function. We assume the solution $\Phi:\mathbb{R}^3 \to \mathbb{R}$ to be smooth over $\mcD$. We will choose $\AAA[\textbf{ }] \equiv \nabla[\textbf{ }] \cdot \nabla[\textbf{ }] = \Delta [\textbf{ }], \BBB[\textbf{ }] \equiv \text{sinh}[\textbf{ }], \CCC \equiv g$, where $g$ is some independent term. Thus, we define:
\begin{equation}
    \label{npbe_defn}
    \FF[\textbf{ }] \equiv \Delta [\textbf{ }] + \text{sinh}[\textbf{ }] + g
\end{equation}
For some solution attempt $\mcN:\mathbb{R}^3 \to \mathbb{R}$, we may write
\begin{equation}
    \label{npbe_attempt}
    \FF(\mcN) = \Delta (\mcN) + \text{sinh}(\mcN) + g 
\end{equation}
Defining $\mcNe = \Phi - \mcN$, and using Taylor expansions
\[\text{sinh}(\Phi) = \text{sinh}(\mcN) + \text{cosh}(\mcN)\mcNe - \text{sinh}(\mcN)\frac{\mcNe^2}{2!} - ...
\]

\[\text{sinh}(\mcN) = \text{sinh}(\Phi) - \left[\text{cosh}(\mcN)\mcNe - \text{sinh}(\mcN)\frac{\mcNe^2}{2!} - ...\right]
\]

\[\text{sinh}(\mcN) = \text{sinh}(\Phi) - \left[\text{cosh}(\mcN)\left(\mcNe + \frac{\mcNe^3}{3!} + ... \right) \right] + \left[\text{sinh}(\mcN)\left(\frac{\mcNe^2}{2!} + \frac{\mcNe ^4}{4!} + ...\right) \right]
\]

\[\text{sinh}(\mcN) = \text{sinh}(\Phi) - \text{cosh}(\mcN)\text{sinh}(\mcNe) + \text{sinh}(\mcN)[1 - \text{cosh}(\mcNe)]
\]
Thus, we may write our originally introduced $\FF(\mcN)$ as:
\begin{equation}
\label{npbe_ec_ver}
\begin{aligned}
&    \FF(\mcN) = \Delta (\Phi) -\Delta (\mcNe) + \text{sinh}(\Phi) - \text{cosh}(\mcN)\text{sinh}(\mcNe) +  \text{sinh}(\mcN)(1 - \text{cosh}(\mcNe)) + g\\
&
\\
& \textbf{ }\textbf{ }\textbf{ }\textbf{ }\textbf{ }\textbf{ } \text{ } = -\Delta (\mcNe) - \text{cosh}(\mcN)\text{sinh}(\mcNe) + \text{sinh}(\mcN)(1 - \text{cosh}(\mcNe))\\
\end{aligned}    
\end{equation}
a transformation from implicit to explicit dependence on $\mcNe$, since $\Delta (\Phi) + \text{sinh}(\Phi) + g = 0$.

Note that we could have used a different method to obtain the same transformation:
\[\text{sinh}(\Phi) = \text{sinh}(\mcN + \mcNe) = \text{sinh}(\mcN)\text{cosh}(\mcNe) + \text{cosh}(\mcN)\text{sinh}(\mcNe)
\]
Adding and subtracting $\text{sinh}(\mcN)$ onto RHS, we get:
\[\text{sinh}(\mcN) = \text{sinh}(\Phi) - \text{cosh}(\mcN)\text{sinh}(\mcNe) + \tsinh(\mcN)[1 - \tcosh(\mcNe)]
\]
which again gives us:
\[\FF(\mcN) = -\Delta (\mcNe) - \text{cosh}(\mcN)\text{sinh}(\mcNe) + \text{sinh}(\mcN)(1 - \text{cosh}(\mcNe))
\]

The same larger goal of obtaining explicit relationships between $\mcL, \mcNe$ is the core pursuit, not the use of a particular kind of transformation.

Now, if we setup a second NN DE solver (while halting the optimization of $\mcN$), we can use Eq. \ref{npbe_ec_ver} to obtain an approximation of $\mcNe$. Just like the original NN DE solver, we have a DE at hand, known constraint conditions, and a NN to solve for the unique solution. Further, $\mcNe \in \mcC^\infty$, which implies it should be in the same \textit{class} of difficulty to approximate as $\Phi$: the error correction NN DE approximation $\hmcNe$ need not consume significantly more resources per iteration, while giving significantly better than expected performance when compared to a converged $\mcN$ estimate for $\Phi$. As such, we not only get an estimate for the error, we also get to use it as a correction term with the same advantages that make $\mcN$ a lucrative option.

\newpage
\section{The Algorithm}
We summarize the general algorithmic approach to internal error correction below:

\begin{algorithm}
\caption{Internal Error estimation and correction}
\begin{enumerate}
    \item Initialize, with random weights and smooth activation, a neural network $\mcN: \mcD \to \mathbb{R}^D$.
    \item Begin optimizing $\mcN$ using: $\mcL = \overline{||\FF[\mcN(x_1)]||^p} + \overline{||\mcNe(x_2)||^p} , \textbf{ }\textbf{ }\textbf{ } x_1 \in \mcD, x_2 \in \pmcD$.
    \item End Step 2, once $\mcL$ performance per iteration consistently drops below a pre-selected standard.
    \item Freeze the parameters of $\mcN$ so that it is now a fixed, well defined, smooth function over $\mcD$.
    \item Use a valid generating transformation for $\BBB'$ (for analytic functions, Taylor expansions are always valid and easily implementable).
    \item Define $\mcNe = \Phi - \mcN$.
    \item Initialize, with random weights and smooth activation, a new neural network $\hmcNe$.
    \item Define the new $\mcL_2 = \overline{||\AAA[\hmcNe(x_1)] + \BBB'[\mcN(x_1), \hmcNe(x_1)] + \FF[\mcN(x_1)]||^{p}} + \overline{||\hmcNe(x_2) - \mcNe(x_2)||^p}$, $x_1 \in \mcD, x_2 \in \pmcD$.
    \item Optimize $\hmcNe$ using $\mcL_2$ until performance drops below some preset standard.
    \item Use $\mcN + \hmcNe$ as an estimation of $\Phi$ over $\mcD$.
\end{enumerate}
\end{algorithm}

\section{Conclusions}
The surging popularity of NN DE solvers presents exciting possibilities in many scientific fields: their capacity to sidestep the curse of dimensionality, the ever decreasing costs of training them, and their easy interpret-ability make for a powerful combo. As such, it is important that the solution estimates they provide be capable of validation over domains where true solutions are not available. Our work proposes one class of methods by which this current deficiency in many types of NN DE solvers may be handled.

Notice that our proposed method does not fully solve the problem: we have simply transferred the ambiguity associated with the estimate $\mcN$ into the ambiguity associated with the estimate $\hmcNe$. However, this is still substantially better than the situation that existed before: even an ambiguous estimate of $\hmcNe$ has clear practical benefits. More importantly, as shown in \cite{dogred20_2, dog20, mar20}, Eq. \ref{rearragned_error_eq} sometimes allows us to obtain formal bounds on $||\mcNe||$ without the need to construct $\hmcNe$: our method can fully solve the ambiguity problem under some extra assumptions. Further, we obtain correction estimates $\hmcNe$ by obtaining Eq. \ref{Err_crrct_gen} from Eq. \ref{gen_loss_DE}. We can repeat the trick to obtain an estimate for $\hmcNe - \mcNe$ by obtaining a new loss $\mcL_3$ just like we obtained $\mcL_2$: indeed we can add corrections to any order under the assumed conditions. However, our results were obtained under very strong assumptions on $\Phi$: future work will focus on DEs where $\Phi$ satisfies weaker conditions.

Lastly, we note that the general ideas we are suggesting are not substantially different than those that already exist for classical numerical methods: higher order corrections are as old as the field itself. By pairing the many significant advantages of modern NNs with those old ideas, we simply hope to have presented a blueprint that will be useful to a wide class of scientific problems that are being tackled using NN DE solvers.

\section{Acknowledgements}
A.S.D. was funded by the President’s PhD Scholarships at Imperial College London. A.S.D.'s research has been supported by the EPSRC Centre for Doctoral Training in Mathematics of Random Systems: Analysis, Modelling and Simulation (EP/S023925/1).

\end{document}